\documentclass[12pt]{article}
 
\usepackage{babel}
\usepackage{isolatin1,times,amssymb,amsmath,latexsym,doublespace}

\usepackage[cmtip,arrow,matrix,curve]{xy}
\CompileMatrices
\newcommand{\pil}{\ar@{-}|(0.496){\object@{>}}|(0.5){\object@{>}}|(0.504){\object@{>}}}
\entrymodifiers={--{\bullet}}
\newcommand{\boej}{12pt}

\newtheorem{theorem}{Theorem}[section]
\newtheorem{lemma}[theorem]{Lemma}
\newtheorem{proposition}[theorem]{Proposition}

\newcommand{\Id}{\mathrm{id}}

\newcommand{\AF}{\mathrm{AF}}

\newtheorem{definition}[theorem]{Definition}
\newtheorem{corollary}[theorem]{Corollary}

\newtheorem{remark1}[theorem]{Remark}
\newenvironment{remark}{\begin{remark1}\rm }{\end{remark1}}

\newcommand{\bproof}{\noindent{\bf Proof: }}
\newcommand{\eproof}{\hfill $\Box$\\}

\newcommand{\w}{\omega}
\renewcommand{\a}{\alpha}
\renewcommand{\b}{\beta}
\newcommand{\vp}{\varphi}
\newcommand{\de}{\delta}

\newcommand{\e}{\varepsilon}
\newcommand{\cP}{{\cal P}}
\newcommand{\cO}{{\cal O}}
\newcommand{\cG}{{\cal G}}

\newcommand{\cK}{{\cal K}}

\newcommand{\bC}{{\mathbb C}}
\newcommand{\bZ}{{\mathbb Z}}
\newcommand{\bN}{{\mathbb N}}

\newcommand{\bT}{{\mathbb T}}
\newcommand{\bR}{{\mathbb R}}

\mathversion{bold} \title{Pure infiniteness, stability and $C^*$-algebras
  of graphs and dynamical systems} 
\mathversion{normal}

\author{Jacob v.B.\ Hjelmborg}

\date{}

\begin{document}
\maketitle

\begin{abstract}
  Pure infiniteness (in sense of \cite{KR:purely}) is considered for
  $C^*$-algebras arising from singly generated dynamical systems. In
  particular, Cuntz-Krieger algebras and their generalizations, i.e.,
  graph-algebras and $\cO_A$ of an infinite matrix $A$, admit
  characterizations of pure infiniteness. As a consequence, these
  generalized Cuntz-Krieger algebras are traceless if and only if they are
  purely infinite.  Also, a characterization of $\AF$-algebras among these
  $C^*$-algebras is given. In the case of graph-algebras of locally finite
  graphs, characterizations of stability are obtained.
\end{abstract}

\section{Introduction}
\label{intro}

Following J.\ Renault \cite{Ren:cuntz-like} a singly generated dynamical
system (SGDS), denoted $(X,T)$, consists in our setting of a locally
compact, second countable Hausdorff space $X$ and a local homeomorphism $T$
from an open subset dom($T$) of $X$ onto an open subset ran($T$) of $X$.
One can associate a locally compact \'etale Hausdorff groupoid, denoted
$G(X,T)$, that is also amenable. The full and reduced $C^*$-algebra of the
convolution algebra of $G(X,T)$ coincide and is named $C^*(X,T)$. These
$C^*$-algebras are of our interest and a goal is to characterize singly
generated dynamical systems $(X,T)$ for which $C^*(X,T)$ is purely infinite
(see definition stated below). This is obtained for the graph-algebras (of
locally finite directed graphs) introduced in \cite{KPRR}, and for the
$C^*$-algebras $\cO_A$, where $A$ is an infinite matrix, introduced by R.\ 
Exel and M.\ Laca in \cite{ExelLaca:infinite}.  The fact that these
$C^*$-algebras come from singly generated dynamical systems is revealed by
J.\ Renault (see \cite[section 4]{Ren:cuntz-like}).

By \cite{KR:purely} a $C^*$-algebra $A$ is {\em purely infinite} if there
are no characters on $A$, and for every pair of positive elements $a,b$ in
$A$ such that $b$ lies in the closed two-sided ideal generated by $a$ it
follows that $b\lesssim a$, i.e., there is a sequence
$(r_j)_{j=1}^{\infty}$ in $A$ such that $r_jbr_j^*\to a$ (cf.\ J.\ Cuntz'
comparison theory for positive elements). An equivalent formulation of pure
infiniteness is contained in \cite{KR:purely}:
\begin{equation}
\label{purely_infinite}
\forall a\in A^+\colon a\oplus a\lesssim a, 
\end{equation}
where $a\oplus a$ denotes the element of $M_2(A)$ with $a$ in the diagonal
and zero elsewhere. The definition agrees with the usual definition by
J.~Cuntz when $A$ is simple, and in the non-simple case the two properties
are neither weaker nor stronger than each other. However, if every non-zero
hereditary $C^*$-subalgebra in every quotient of $A$ contains an infinite
projection, then $A$ is purely infinite, as shown in \cite{KR:purely}. It
is an open question whether a $C^*$-algebra $A$ with no non-zero trace
(defined on a not necessarily closed or dense two-sided ideal in $A$) is
purely infinite. The answer is affirmative for the $C^*$-algebras
considered in this paper.

In section \ref{purely} we consider locally finite directed graphs which
may have infinitely many vertices, but only finitely many edges in and out
of each vertex. In \cite{KPRR} it is shown that to each such graph $E$ one
can associate a locally compact groupoid $\cG_E$ such that its groupoid
$C^*$-algebra $C^*(\cG_E)$ is the universal $C^*$-algebra generated by
families of partial isometries satisfying Cuntz-Krieger relations
determined by $E$. In fact this $C^*$-algebra comes from an SGDS
$(X_{A_E},T_{A_E})$, where $T_{A_E}$ is the one-sided shift on the infinite
path space $X_{A_E}$ associated to the edge-matrix $A_E$ of the graph $E$
(cf.\ \cite{KPRR} and section \ref{sgds}). Pure infiniteness is
characterized for these $C^*$-algebras in Theorem \ref{thm-1}. A
$C^*$-algebra $A$ is stable if $A\cong A\otimes\cK$, where $\cK$ is the
$C^*$-algebra of compact operators on a separable infinite dimensional
Hilbert space. Theorem \ref{main_theorem} give characterizations of
stability for $C^*$-algebras of locally finite directed graphs. In
particular, stability is characterized by the absence of (non-zero) unital
quotients and bounded trace. Using this characterization it is shown that
the extension of stable $C^*$-algebras of graphs is stable (cf.\ Corollary
\ref{thm3-1}). Also, a graph-theoretical analogue of stability is given.

In section \ref{groupoids} a sufficient condition on $r$-discrete
essentially principal groupoids is given for which the associated reduced
$C^*$-algebra is purely infinite, and it is shown in section
\ref{pure_infiniteness} that this condition is also necessary in case of
the generalized Cuntz-Krieger algebras associated to infinite matrices
which are not necessarily row-finite (cf.\ Theorem \ref{markov_pi}). In
section \ref{infinite_matrices}, the $\AF$-algebras among these
$C^*$-algebras are characterized (cf.\ Theorem \ref{af}).

The author is very grateful for helpful discussions with M.\ R{\o}rdam, A.\ 
Kumjian, C.\ Anantharaman-Delaroche and J.\ Renault and would like to thank
the operator algebras group (or groupoid ?) in Orleans for kind hospitality
during a visit in March '99.

\mathversion{bold}
\section{Pure infiniteness and stability of $C^*$-algebras of locally
  finite graphs.} 
\mathversion{normal}
\label{purely}
\subsection{Preliminaries}
\label{preliminaries}
A {\em directed graph} $E$ is a quadruple $(E^0,E^1,r,s)$ where $E^0$
and $E^1$ are countable sets and $r,s\colon E^1\to E^0$ are maps. $E^0$
is the set of vertices and $E^1$ is the set of edges in the
graph. The maps $s$ and $r$ are called the source and range maps and
they give the direction of each edge. The graph $E$ is {\em row-finite}
if for every $v\in E^0$ the set $s^{-1}(v)\subseteq E^1$ is finite,
i.e., only finitely many edges come out of each vertex, and if in
addition $r^{-1}(v)$ is finite for all $v\in E^0$ then $E$ is {\em
locally finite}. A finite \emph{path} in $E$ is a sequence $\a = 
(e_1,e_2,\ldots,e_k)$ of edges in $E$ with $s(e_{j+1})=r(e_j)$ for 
$1\le j \le k-1$. We write $E^*$ for the finite path space of $E$, 
and $E^{\infty}$ for the infinite path space.
The range and source maps $r,s$ extend naturally to $E^*$.
Suppose $\a$ is a path in $E$. The subset $\a^0\subseteq E^0$ given by
\[
\a^0 =\{v\in E^0\colon\mbox{$v=s(e)$ or $v=r(e)$ for some $e\in \a$}\},
\]
is the set of vertices of $\a$.  A {\em detour} of a path $\a$ is a path
$\beta \in E^*$ such that $s(\beta)\in \a^0$ and $r(\beta)\in \a^0$.
Trivially, subpaths of $\a$ are detours.  A path $\a\in E^*$ of lenght
$|\a|>0$ is a cycle based at $v$ if $s(\a)=v=r(\a)$. A vertex $v\in E^0$
that emits no edges is called a \emph{sink}. We now consider sets of
vertices that give rise to ideals in the associated $C^*$-algebra (consult
\cite{KPRR} for details). Let $H$ be a subset of $E^0$. $H$ is {\em
  hereditary} if, whenever $v\in H$ and there is a path from $v$ to $w$,
then $w\in H$.  $H$ is {\em saturated} if $v\in H$ when every edge $e\in
E^1$ with source $s(e)=v$ has range $r(e)\in H$.

Given a row-finite directed graph $E$, a {\em Cuntz-Krieger $E$-family}
is given by a set of mutually orthogonal projections indexed by
vertices, $\{P_v\colon v\in E^0\}$ and a family of partial isometries
indexed by edges $\{S_e\colon e\in E^1\}$ subject to the following
conditions:
\begin{eqnarray}
S^*_eS_e = P_{r(e)}
\label{graph1}
\end{eqnarray}
for all $e\in E^1$ and
\begin{eqnarray}
P_v = \sum_{\{e\colon s(e)=v\}}S_eS_e^*
\label{graph2}
\end{eqnarray}
for all $v\in s(E^1)$. The universal $C^*$-algebra generated by a
Cuntz-Krieger $E$-family $\{ P_v,S_e\}$ is denoted by $C^*(E)$ (See 
\cite[Theorem 1.2]{KPR}). Note that the projections $\{P_v\}$ are the 
initial projections of the partial isometries $S_e$ with $r(e)=v$, and
note also that a projection $P_v$ can be non-zero even if there are no 
edges coming out of $v$.

Recall that a projection $p$ in a $C^*$-algebra $A$ is {\em infinite} if it
is equivalent to a proper subprojection of itself, and $p$ is {\em finite}
otherwise. If there exist mutually orthogonal subprojections $q$ and $r$ of
$p$ such that $q$ and $r$ are both equivalent to $p$, then $p$ is said to be
{\em properly infinite}. Equivalently, $p$ is properly infinite if $p\oplus
p \lesssim p$ (and $p$ is non-zero). 

A trace on $A$ is a positive linear map $\tau\colon I\to\bC$, where
$I$ is a (not necessarily closed or dense) two-sided ideal in $A$, such that
$\tau(x^*x)=\tau(xx^*)$ for all $x\in I$. Notice that if
$\tau$ is a bounded trace on $I$ then it extends to a bounded trace
on $A$.

\mathversion{bold}
\subsection{Stable $C^*$-algebras of locally finite graphs.} 
\mathversion{normal}
\subsubsection{Stability and projections}
\label{section charac_proj}

By \cite[Theorem3.3]{HjeRor:stable} a $C^*$-algebra $A$ that admits a
countable approximate unit consisting of projections is stable if and only
if for each projection $p\in A$ there is a projection $q\in A$ such that
$p\sim q$ and $p\perp q$. Lemma \ref{lemma_charac} below is a reformulation
suited for $C^*$-algebras of graphs.

\begin{lemma}
\label{lemma_app-unit}
Let $A$ be a $C^*$-algebra equipped with an increasing approximate
unit $(p_n)^\infty_{n=1}$ consisting of projections. The
following statements are equivalent:
\begin{itemize}
\item[(i)] $\forall\, n\in\bN \,\exists\, m>n\,\colon p_n\lesssim p_m-p_n$.
\item[(ii)] $\forall\, p\in P(A)\,\exists\, q\in P(A)\,\colon p\sim q$ and
$p\perp q$.
\end{itemize}
\end{lemma}

\bproof

{\bf (i) $\Rightarrow$ (ii):}  Let $p\in A$ be a projection. By
  \cite[Lemma 3.1]{HjeRor:stable} there is $k\in\bN$, $\tilde{p}_k\in P(A)$
and a unitary $u\in\tilde{A}$ such that
\[
p\le\tilde{p}_k,\quad \tilde{p}_k = u^*p_ku
\]
From this and the assumption it follows that
\[
upu^*\le u\tilde{p}_ku^* = p_k\lesssim p_m-p_k
\]
for some $m>k$. Hence there is a projection $\tilde{q}\in P(A)$ such that
\[
upu^*\sim\tilde{q},\quad\tilde{q}\perp upu^*
\]
Put $q=u^*\tilde{q}u$ and we are done.

{\bf (ii) $\Rightarrow$ (i):} Let $n\in\bN$. 
By assumption there is $q\in P(A)$ satisfying
$q\sim p_n$ and $q\perp p_n$. Note that $q\in(1-p_n)A(1-p_n)$. For
each $k\ge n$ set $\hat{p}_k = p_k-p_n$. Then $(\hat{p}_k)^\infty_{k=n}$
is an increasing approximate unit for $(1-p_n)A(1-p_n)$. By
\cite[Lemma 3.1]{HjeRor:stable} there is $m\in\bN$ and a projection $r \in
(1-p_n)A(1-p_n)$ such that
\[
q\le r,\quad r\sim\hat{p}_m
\]
Hence
\[
p_n\sim q\lesssim\hat{p}_m = p_m-p_n
\]
and the proof is complete.
\eproof

\begin{remark} 
Let $E$ be a directed graph and let $\{P_v,S_e\}$ be the
associated Cuntz-Krieger $E$-family. The set of vertices $E^0$ is infinite
if and only if $C^*(E)$ is non-unital (\cite[Proposition 1.4]{KPR}).
Suppose that $E^0=\{v_j\colon j\in\bN\}$ is infinite and
define $p_n = \sum^n_{j=1} P_{v_j}$ for $n\in\bN$. Then
$(p_n)^\infty_{n=1}$ is a strictly increasing approximate unit for
$C^*(E)$ consisting of projections. 
\end{remark}

\begin{lemma}
\label{lemma_charac}
Let $E$ be a row-finite directed graph. Then $C^*(E)$ is stable if
and only if for each finite subset $F\subseteq E^0$ there is a finite
subset $G\subseteq E^0$ such that $F\cap G=\emptyset$ and $\sum_{v\in F}
P_v\lesssim\sum_{v\in G}P_v$. 
\end{lemma}

\bproof
If $C^*(E)$ is stable then it is non-unital, hence $E^0$ is infinite and
$(p_n)^\infty_{n=1}$ defined above is a strictly increasing approximate
unit  for $C^*(E)$ consisting of projections. By \cite[Theorem
3.3]{HjeRor:stable} and Lemma \ref{lemma_app-unit} the desired property
holds. Conversely, if for each finite subset $F\subseteq E^0$ there is a
finite subset $G\subseteq E^0$ such that $F\cap G=\emptyset$ and $\sum_{v\in F}
P_v\lesssim\sum_{v\in G}P_v$ then $(p_n)^\infty_{n=1}$ satisfies statement
(i) of Lemma \ref{lemma_app-unit} and $C^*(E)$ is stable by \cite[Theorem
3.3]{HjeRor:stable}.
\eproof

\subsubsection{Lifting}
\label{subsection_lift}

This section deals with the lifting of equivalent projections in a
quotient of a $C^*$-algebra (cf.\ Lemma~\ref{lemma_lift}). See
\cite[Lemma 4.11]{KR:purely} for a proof of the following Lemma.

\begin{lemma}
\label{lemma_exact}
Let $0\longrightarrow I\longrightarrow A\stackrel{\pi}{\longrightarrow}
B\longrightarrow 0$ be an extension of $C^*$-algebras and let $e,f$
be projections in $A$. Then $\pi(e)\lesssim\pi(f)$ if and
only if there exists $x\in I^+$ such that $e\lesssim f\oplus x$.
\end{lemma}

In the following $A$ is a $C^*$-algebra, $\cP=\{q_n\}_{n=1}^\infty$ is
an increasing sequence of projections in $A$, and $I(\cP)$ denotes the
ideal in $A$ generated by $\cP$.

\begin{lemma}
\label{lemma_ideal}
$I(\cP)\otimes\cK$ has an approximate unit consisting of projections.
\end{lemma}

\bproof
Put $B = \overline{\bigcup^\infty_{n=1} q_nAq_n}$. Then $B$ is a full
hereditary subalgebra of $I(\cP)$. By Brown's Theorem
\cite{Bro:stabher}, $B\otimes\cK\cong I(\cP)\otimes\cK$ and
$B\otimes\cK$ clearly admits an approximate unit consisting of projections.
\eproof

\begin{lemma}
\label{lemma_lift}
Assume $\cP$ satisfies
\begin{equation}
\label{eq*}
\forall p\in\cP\,\exists q\in\cP\colon q\ge p\ \mbox{and} \ p\lesssim q-p.
\end{equation}
Let $\pi\colon A\to A/I(\cP)$ be the quotient homomorphism. Let $e,f\in
A$ be projections such that $\pi(e)\lesssim \pi(f)$. Then
$e\lesssim f\oplus p$ for some projection $p\in\cP$.
\end{lemma}

\bproof
By Lemma \ref{lemma_exact} there exists $x\in I(\cP)^+$ such that
$e\lesssim f\oplus x$. 
Embed $I(\cP)$ into $I(\cP)\otimes\cK$. For each $\e>0$ define
$\vp_\e\colon\bR_+\to\bR_+$ by
\[
\vp_\e(t) = \left\{\begin{array}{lll} 0 & \mbox{for} & t\le\e,\\
\e^{-1}(t-\e) & \mbox{for} & \e\le t\le 2\e,\\ 1 & \mbox{for} & t\ge
2\e.\end{array}\right. 
\]
By \cite[Prop. 2.4]{Ror:UHFII} there is $\de\in{}]0,\frac12[$ such that
\[
e\lesssim f\oplus \vp_{\de}(x).
\]
By Lemma \ref{lemma_ideal} there is a projection $q\in I(\cP)\otimes\cK$
such that
$\|x-qxq\|<\de$. By \cite[Prop. 2.2]{Ror:UHFII}, which generalizes the
fact that close projections are equivalent, $\vp_{\de}(x)\lesssim
qxq$. Hence
\[
e\lesssim f\oplus\vp_{\de}(x)\lesssim f\oplus qxq\lesssim f\oplus q.
\]
From the description of $I(\cP)\otimes\cK$ in Lemma \ref{lemma_ideal}
there exists
$p'\in\cP$ and $m\in\bN$ such that $q\lesssim p'\otimes 1_m$. By the
assumed property, (\ref{eq*}), there exists $p\in\cP$ such that
$p'\otimes 1_m\lesssim p$. This proves that $e\lesssim f\oplus p$.
\eproof

\subsubsection{Graph-trace}
\label{subsection_graph-trace}

\begin{definition}
\label{definition graph-trace}
Let $E$ be a directed graph. A map $\tau_E\colon E^0\to\bR^+$ 
defined on the set of vertices of $E$ is said to be a
graph-trace if
\[
\tau_E(v) = \sum_{\{e\colon s(e)=v\}} \tau_E(r(e)),
\]
and it is bounded if
\[
\sum_{v\in E^0} \tau_E(v) < \infty.
\]
\end{definition}

\begin{lemma}
\label{lemma_graph-trace}
Let $E$ be a locally finite directed graph with no cycles. Then every
bounded graph-trace on $E$ comes from a bounded trace on $C^*(E)$.
\end{lemma}

\bproof Let $\tau_E$ be a bounded graph-trace on $E$. Since $E$ has no
cycles $C^*(E)$ is $\AF$ by Theorem \ref{af} (see Remark \ref{af_remark}).
Furthermore there is an increasing sequence $E_1\subseteq
E_2\subseteq\cdots$ of finite subgraphs of $E$ such that
\[
C^*(E) = \varinjlim(C^*(E_n),\iota_n)
\]
where $\iota_n$ are inclusion maps (see \cite[Remark 2.5]{KPR}). For each
$n\in\bN$, $C^*(E_n)$ is finite dimensional and if $v_1,\dots,v_{k_n}$
denote the sinks of $E_n$ (there is at least one sink) and $m(v_i)$ is
the number of paths with range $v_i$ for $i=1,\dots,k_n$ then since $E$ 
has no cycles it follows from \cite[Corollary 2.3]{KPR} that
\begin{eqnarray}
\label{description}
C^*(E_n) &\cong & \bigoplus^{k_n}_{i=1} M_{n(v_i)}(\bC),
\end{eqnarray}
where the isomorphism in (\ref{description}) maps the projection
$P_{v_i}\in C^*(E_n)$ to a minimal projection in $M_{n(v_i)}(\bC)$.
Hence for each $n\in \bN$ there is a unique trace
$\tau_n\colon C^*(E_n)\to\bC$ such that
$$
\tau_n(P_{v_i})=\tau_E(v_i),\qquad i=1,\ldots,k_n
$$
Also from \cite[Corollary 2.3]{KPR}
\begin{eqnarray*}
K_0(C^*(E_n)) &=& \langle\{[P_{v_i}]\colon i=1,\dots,k_n\}\rangle\cong
\bZ^{k_n},
\end{eqnarray*}
and for every $v\in E^0$
\begin{eqnarray}
\label{proj}
[P_v]=\sum_{i=1}^{k_n}n_{v_i}[P_{v_i}] 
\end{eqnarray}
where $n_{v_i}$ is the number of paths from $v$ to $v_i$. Let $v\in
E_n^0$. Since $E_n$ is finite with no cycles an induction using the
graph-trace relation $\tau_E(v) = \sum_{\{e\colon s(e)=v\}} \tau_E(r(e))$
reveals that
\begin{eqnarray}
\label{trace}
\tau_E(v)=\sum_{i=1}^{k_n}n_{v_i}\tau_E(v_i). 
\end{eqnarray}
Hence from (\ref{proj}) and (\ref{trace}) we get that 
$$
\tau_n(P_v)=\sum_{i=1}^{k_n}n_{v_i}\tau_n(P_{v_i})=
\sum_{i=1}^{k_n}n_{v_i}\tau_E(v_i)=\tau_E(v),
$$
for all $v\in E_n^0$
By uniqueness $\tau_{n+1}$ extends $\tau_{n}$. Define $\tau$ on the
$*$-subalgebra $\bigcup^\infty_{n=1} C^*(E_n)$ of $C^*(E)$ by setting
$\tau(a)=\tau_n(a)$ if $a\in C^*(E_n)$. Put $k=\sum_{v\in E^0}\tau_E(v)$. If
$a\in \bigcup^\infty_{n=1} C^*(E_n)$ then 
\begin{eqnarray*}
|\tau(a)| &\le& \sum_{v\in E^0}\tau(P_v)\|a\| \\
&=& \sum_{v\in E^0}\tau_E(v)\|a\| \\
&=& k \|a\|.
\end{eqnarray*}
Hence $\|\tau\|< \infty$ and $\tau$ is a bounded linear
function from $\bigcup^\infty_{n=1} C^*(E_n)$ to $\bC$, so it extends to
a bounded linear functional $\tilde{\tau}$ on $C^*(E)$. It is easily
checked that $\tilde{\tau}$ is a trace on $C^*(E)$. This completes the proof.
\eproof

\begin{remark}
Obviously the converse statement of Lemma \ref{lemma_graph-trace}
holds, i.e., a trace on $C^*(E)$ induces a graph-trace on $E$. The
conclusion of Lemma \ref{lemma_graph-trace} is likely to hold more
general than for directed graphs with no cycles.
\end{remark}

\begin{lemma}
\label{proposition no cycles}
Let $E$ be a locally finite directed graph with no cycles. Then $C^*(E)$
is stable if and only if $E$ has no non-zero bounded graph-trace.
\end{lemma}

\bproof Since $E$ has no cycles, $C^*(E)$ is $\AF$ by Theorem \ref{af} (Remark
\ref{af_remark}). Hence $C^*(E)$ is stable if and only if $C^*(E)$ admits no
non-zero bounded trace by \cite{Bla:traces}. By Lemma
\ref{lemma_graph-trace}, $C^*(E)$ admits no non-zero bounded trace if and
only if $E$ admits no non-zero bounded graph-trace.  \eproof

\subsubsection{Stability}
\label{subsection_main}

\begin{definition}
\label{left-infinite vertex}
Let $E=(E^0,E^1,r,s)$ be a locally finite directed graph. A vertex $v\in
E^0$ is left-infinite if the set
\[
L(v) = \{w\in E^0\colon \mbox{$s(\alpha)=w$ and $r(\alpha)=v$ for some
  $\alpha\in E^*$}\}
\]
is infinite. $v$ is left-finite if $L(v)$ is finite.
A subgraph $F\subseteq E$ is left-infinite (resp.\ left-finite) if every
vertex in $F^0$ is left-infinite (resp.\ left-finite).
\end{definition}

\begin{lemma}
\label{konig}
Let $E$ be a locally finite directed graph. A vertex $v\in E^0$ is
left-infinite if and only if there is an infinite path $\a\in E^{\infty}$
such that all edges of $\a$ are distinct and $r(\a)=v$.
\end{lemma}

\bproof 

{\bf only if:} This may be recognized by graph-theorists as K{\"o}nig's infinity
lemma \cite{Konig}. Let $v\in E^0$ be left-infinite. Since $L(v)$ is
infinite  and only
finitely many edges enter $v$ there is an edge $e_1\in E^1$ which is the
last edge in infinitely many paths $\alpha\in E^*$ with $s(\alpha)\in L(v)$
and $r(\alpha)=v$. Furthermore we may assume that $s(e_1)\neq r(e_1)$. Now,
since only finitely many edges enter $s(e_1)$ there is an edge $e_2\in E^1$
with $r(e_2)=s(e_1)$, $s(e_2)\neq s(e_1)$, $s(e_2)\neq v$ and such that the
path $(e_2,e_1)$ can be concatenated to infinitely many paths $\alpha\in
E^*$ with $s(\alpha)\in L(v)$ and $r(\alpha)=s(e_2)$. An obvious induction
gives an infinite sequence $(e_i)_{i=1}^\infty$ of edges in $E^1$
with the desired properties.

{\bf if:} This follows easily.
\eproof

\begin{lemma}
\label{lemma_left-infinite}
Let $E$ be a locally finite directed graph. If $E$ is left-infinite then
$C^*(E)$ is stable.
\end{lemma}

\bproof
Note that for every edge $e\in E^1$,
\[
P_{r(e)}=S^*_{e}S_{e}\sim S_{e}S^*_{e}\le P_{s(e)},
\]
and $P_{r(e)}\perp P_{s(e)}$ if $r(e)\neq s(e)$.
By Lemma~\ref{konig} there is for every vertex $v\in E^0$ a sequence
$\{e_i\}_{i=1}^\infty$ of edges such that the vertices $\{r(e_i)\}_{i=1}^\infty$
are mutually distinct, $r(e_{i+1})=s(e_i)$ and $r(e_1)=v$. 
Hence there is a sequence
$(P_{s(e_i)})_{i=1}^\infty$ of mutually orthogonal projections satisfying
\[
P_v\lesssim P_{s(e_1)}\lesssim P_{s(e_2)}\cdots ,
\]
for all $v\in E^0$. This implies that for each finite subset
$F\subseteq E^0$ there is a finite 
subset $G\subseteq E^0$ such that $F\cap G=\emptyset$ and $\sum_{v\in F}
P_v\lesssim\sum_{v\in G}P_v$, hence $C^*(E)$ is stable by
Lemma~\ref{lemma_charac}. 
\eproof

\begin{theorem}
\label{main_theorem}
Let $E$ be a locally finite directed graph with no sinks. The following
five statements are equivalent
\begin{itemize}
\item[(a)] $C^*(E)$ is stable.
\item[(b)] $C^*(E)$ admits no non-zero unital quotient and no non-zero
bounded trace.
\item[(c)] No non-zero quotient of $C^*(E)$ is either unital or an
$\AF$-algebra with a non-zero bounded trace.
\item[(d)] $E$ has no left-finite cycles and the subset $S^0\subseteq E^0$
given by  
\[
S^0=\{v\colon \mbox{$v\in {\gamma}^0$ and ${\gamma}^0$ is left-finite for
  some infinite path $\gamma\in E^\infty$}\}
\] 
admits no non-zero bounded graph-trace (cf.\ section
\ref{subsection_graph-trace}).
\item[(e)] For every finite subset $V\subseteq E^0$ there exists a    
finite subset $W\subseteq E^0$ disjoint to $V$ such that 
\[
\sum_{v\in V}P_v\lesssim \sum_{v\in W}P_v. 
\]
\end{itemize} 
\end{theorem}

\bproof

{\bf (a) $\Rightarrow$ (b):} A stable $C^*$-algebra is non-unital with no
non-zero bounded trace.  Every quotient of a stable $C^*$-algebra is
stable, hence non-unital (cf.\ Corollary \ref{thm3-1}).

{\bf (b) $\Rightarrow$ (c):} No non-zero quotient of $C^*(E)$ admits
  a bounded trace since it would lift.

{\bf (c) $\Rightarrow$ (d):} First we  show that if $E$ has a
left-finite cycle, then $C^*(E)$ has a unital quotient. Consider a
left-finite cycle $\a$, i.e. a cycle with left-finite vertices, and let $v\in
{\a}^0$. Define a subset $L(v)\subseteq E^0$ by
\[
L(v) = \{w\in E^0\mid \mbox{$s(\a)=w$ and $r(\a)=v$ for some
  $\a\in E^*$}\},
\]
i.e., $L(v)$ contains all vertices $w$ for which there is a path connecting
$\w$ to $v$. It follows that $L(v)$ is finite since $v$ is left-finite.
Set $H^0=E^0\backslash L(v)$. If $H^0=\emptyset$ then $E^0$ must be finite,
hence $C^*(E)$ is unital.  Suppose that $H^0$ is not empty. $H^0$ is a
hereditary and saturated subset of $E^0$; $H^0$ is hereditary because if
not we may assume there exist $u\in H^0$, $w\in L(v)$ and an edge with
source $u$ and range $w$. But then there is a path connecting $u$ to $v$,
hence $u\in L(v)$ and this is a contradiction. $H^0$ is saturated because
for every $w\in L(v)$ there is an edge $e\in E^1$ such that $r(e)\in L(v)$
and $s(e)=w$, since $v\in {\a}^0$. Since $E$ has no sinks and $H^0$ is a
hereditary and saturated subset of $E^0$, there is by \cite[Theorem
6.6]{KPRR} a two-sided closed ideal $I(H^0)$ of $C^*(E)$ such that the
quotient $C^*$-algebra $C^*(E)/I(H^0)$ is naturally isomorphic to $C^*(F)$
of the directed graph $F= (L(v),\{e\in E^1\colon r(e)\in L(v)\})$. $C^*(F)$
is unital since $L(v)$ is finite, hence $C^*(E)$ has a unital quotient.

Consider the subset $S^0$ defined in item (d).  We show that $S^0$ induces
an $\AF$-quotient of $C^*(E)$.  Set $H^0 = E^0\backslash S^0$.  $H^0$ is a
hereditary and saturated subset of $E^0$: Suppose $H^0$ is not hereditary,
i.e. there exist $v\in H^0$, $w\in S^0$ and an edge with source $v$ and
range $w$. $v$ is not left-finite because otherwise $v\in S^0$, which is a
contradiction. But then $v$ must be left-infinite, which implies that $w$
is left-infinite. This is also a contradiction, so $H^0$ must be
hereditary. $H^0$ is saturated because if $w\in S^0$ then there is an edge
$e\in E^1$ such that $r(e)\in S^0$ and $s(e)=w$.  Since $E$ has no sinks
and $H^0$ is hereditary and saturated we obtain by \cite[Theorem 6.6]{KPRR}
that
\[
C^*(E)/I(H^0)\cong C^*(S),
\]
where $S$ is the subgraph of $E$ given by $S=(S^0,\{e\in E^1\colon
r(e)\notin H^0\})$.  As shown above $E$ has no left-finite cycles, hence
$S$ has no cycles.  Therefore $C^*(S)$ is an $\AF$-algebra by Theorem
\ref{af} (Remark \ref{af_remark}), and by Lemma \ref{lemma_graph-trace}
every bounded graph-trace on $S$ comes from a bounded trace on $C^*(S)$. By
assumption no such exist.  This proves statement (d).

{\bf (d) $\Rightarrow$ (e):} Suppose $S^0\ne \emptyset$ and put $H^0=
E^0\backslash S^0$. $H^0$ is a hereditary and saturated subset of $E^0$ by
an argument similar to that given in the proof of statement (c) implies
(d). We make the following
\begin{description}
\item[Claim:] For every finite subset $K\subseteq H^0$ there is a finite
  subset $M\subseteq H^0$ such that $K\cap M=\emptyset$ and $\sum_{v\in
    K}P_v\lesssim\sum_{v\in M}P_v$.
\end{description}
{\bf Proof of claim:}
It is enough to show that for every vertex $v\in H^0$ there
are left-infinite vertices $w_1,\dots,w_n$ in $H^0$ for some $n\in\bN$ such
that
\[
P_v \lesssim P_{w_1}+\cdots + P_{w_n}
\]
since for each $j=1,\dots,n$ there is an infinite subset
$\{w_j^{(i)}\colon i\in\bN\}$ of $H^0\backslash\{w_j\}$ such that
\[
P_{w_j}\lesssim P_{w_j^{(1)}}\lesssim P_{w_j^{(2)}}\lesssim\cdots
\]
as shown in the proof of Lemma \ref{lemma_left-infinite}.

If $v\in H^0$ is left-infinite we are done. Suppose $v\in H^0$ is
left-finite and 
\[
P_v = \sum_{e\in G_v} S_eS_e^*,\quad G_v = \{e\in E^1\colon s(e)=v\}.
\]
Note that $G_v$ is finite since $E$ is row-finite. For edges $e\in E^1$ let
$m(r(e))$ denote the number of edges with source $v$ and range $r(e)$.
Since $P_{r(e)}=S^*_eS_e$ it follows that
\begin{equation}
\label{equivalent}
P_v\sim \sum_{\{r(e)\colon e\in G_v\}} P_{r(e)}\otimes 1_{m(r(e))}.
\end{equation}
Consider the vertices $r(e)$ for $e\in G_v$.  If $r(e)$ is left-infinite
then there is a subset $\{v_i\colon i\in\bN\}\subseteq~H^0$ of
left-infinite vertices such that
\[
P_{r(e)}\lesssim P_{v_1}\lesssim P_{v_2}\lesssim\cdots
\]
It follows that there is a sequence $(F_j)_{j\in\bN}$ of finite mutually
disjoint subsets of $H^0$ such that
\begin{equation}
\label{majorize}
P_{r(e)}\otimes 1_{m(r(e))}\lesssim\sum_{w\in F_j} P_w
\end{equation}
for each $j\in\bN$. If every $r(e)$ is left-infinite for $e\in G_v$ then
(\ref{equivalent}) and (\ref{majorize}) gives the desired property.  If
$r(e)$ is not left-infinite then we repeat the argument with $r(e)$ in
place of $v$. Recall that $E$ has no left-finite cycles nor sinks.  We
continue repeating the argument until all vertices involved are
left-infinite. The process stops since otherwise $v\in S^0$ which is a
contradiction.  It is now straightforward to show that there is a sum of
projections, associated to left-infinite vertices of $H^0$, that contain a
subprojection equivalent to $P_v$. This proves the claim.\\

Notice that $H^0$ is either an infinite set or the empty set. Suppose
$H^0=\{v_j\colon j\in\bN\}$, $p_n=\sum^n_{j=1} P_{v_j}$ for $n\in\bN$ and
set $\cP=\{p_n\}^\infty_{n=1}$. Let $I(H^0)=\mbox{ideal}(\cP)$ be the
two-sided closed ideal of $C^*(E)$ generated by $\cP$. $H^0$ is hereditary
and saturated and since $E$ has no sinks it follows by \cite[Theorem
6.6]{KPRR}, that the quotient $C^*$-algebra $C^*(E)/I(H^0)$ is isomorphic
to $C^*(S)$, where $S=(S^0,\{e\in E^1\colon r(e)\notin H^0\})$.  Since $S$
has no left-finite cycles and $S^0$ admits no non-zero bounded graph-trace,
$C^*(S)$ is stable by Lemma \ref{proposition no cycles}.

We have shown that if $S$ exists then $C^*(S)$ is stable and isomorphic to
a quotient of $C^*(E)$. If $S=\emptyset$ then $E^0=H^0$ and the claim above
gives that statement (e) holds. If $H^0$ is the empty set (this occurs when
$E$ has no left-infinite vertices) then $C^*(E)$ equals $C^*(S)$ hence is
stable and statement (e) follows by Lemma \ref{lemma_charac}.

Let $V\subseteq E^0$ be a finite subset of $E^0$ and consider the
decomposition
\begin{equation}
\label{decomposition}
\sum_{v\in V}P_v = \sum_{v\in V\cap S^0} P_v + \sum_{v\in V\cap H^0} P_v.
\end{equation}
Since $C^*(S)$ is stable there is a finite subset $W\subseteq S^0$ such
that $W\cap V=\emptyset$ and
\[
\pi\bigg(\sum_{v\in V\cap S^0} P_v\bigg)\lesssim\pi\bigg(\sum_{v\in
W}P_v\bigg),
\] 
where $\pi\colon C^*(E)\to C^*(S)$ is the homomorphism induced
by the quotient map $C^*(E)\to C^*(E)/I(H^0)$. With $\cP$ given above
$C^*(E)$ satisfies property (\ref{eq*}) of Lemma
\ref{lemma_lift} by the claim. Hence by Lemma
\ref{lemma_lift} there is a finite subset $M\subseteq H^0$ such that
\begin{equation}
\label{bridge}
\sum_{v\in V\cap S^0} P_v\lesssim \sum_{v\in W} P_v + \sum_{v\in M} P_v.
\end{equation}
By the claim there is a finite subset $K\subseteq H^0$ such that
$K\cap (V\cup M)=\emptyset$ and
\begin{equation}
\label{from_claim}
\sum_{v\in H^0\cap(V\cup M)} P_v\lesssim \sum_{v\in K} P_v.
\end{equation}
Combining (\ref{decomposition}), (\ref{bridge}) and (\ref{from_claim})
we obtain that the subsets $V$, $W$ and $K$ are mutually disjoint and
\[
\sum_{v\in V} P_v\lesssim \sum_{v\in W} P_v+\sum_{v\in K}P_v.
\]
This proves statement (e).

{\bf (e) $\Rightarrow$ (a):} This follows from Lemma \ref{lemma_charac} and
the proof is complete.
\eproof

\begin{corollary}
\label{thm3-1}
Given an extension
\[
0 \longrightarrow I \stackrel{\displaystyle{\varphi}}{\longrightarrow} A
\stackrel{\displaystyle{\pi}}{\longrightarrow} B \longrightarrow 0,
\]
where $I$, $A$ and $B$ are $C^*$-algebras of locally finite graphs without
sinks. Then $A$ is stable if and only if $I$ and $B$ are stable.
\end{corollary}

\bproof
Assume that $A$ is stable and that $I$ is a (two-sided closed) ideal of
$A\otimes\cK$. The ideal $I$ equals $J\otimes\cK$ for some ideal $J$ of
$A$. This can be seen as follows; Note that $A\otimes\cK$ is $^*$-isomorphic to
the inductive limit $C^*$-algebra $(\bigcup M_n(A))^-$, where $M_n(A)$ is
embedded into the upper left-hand corner of $M_{n+1}(A)$. A computation
reveals that every ideal of $M_n(A)$ is of the form $M_n(J)$
for some ideal $J\subseteq A$ for every $n\in\bN$. Hence every ideal of 
$(\bigcup M_n(A))^-$ is of the form $(\bigcup M_n(J))^-$ for some ideal
$J\subseteq A$.

Let $\pi\colon A\to A/J$ be the quotient $^*$-homomorphism. The
$^*$-homomorphism $\pi\otimes\Id$ of $A\otimes\cK$ into
$(A/J)\otimes\cK$ is surjective with $J\otimes\cK$ as kernel. It follows that
$B\cong A\otimes\cK/I\cong (A/J)\otimes\cK$, hence $B$ is stable.

Assume now that $I$ and $B$ are stable. Then $I$ and $B$ admit no
non-zero bounded trace and no non-zero unital quotient by Theorem
\ref{main_theorem}. $A$ admits no 
non-zero bounded trace either since if $\tau$ is such a trace, then $\tau$
must be zero on $I$, hence $\tau$ induces a non-zero bounded trace on $B$.
Suppose to reach a contradiction that $A$ contains a two-sided closed ideal
$J$ such that the quotient $A/J$ is non-zero and unital. There is a
commuting diagram:
\[
\begin{array}{ccccccccc}
0 & \longrightarrow & I & \longrightarrow & A &
\longrightarrow & A/I & \longrightarrow & 0 \\
&& \downarrow \pi && \downarrow \varphi && \downarrow \psi \\
0 & \longrightarrow & I/(I\cap J) & \longrightarrow & A/J &
\longrightarrow & A/(I+J) & \longrightarrow & 0, 
\end{array}
\]
where $\pi$, $\varphi$ and $\psi$ are surjective $^*$-homomorphisms, and
where each row is exact. It follows that the quotient $A/(I+J)$ is unital
and if $A/(I+J)$ is zero then $I/(I\cap J)$ is $^*$-isomorphic to $A/J$, hence 
non-zero and unital.
This contradicts that $I$ and $B$ admit no non-zero unital quotients. It
follows from Theorem \ref{main_theorem} that $A$ is stable. 
\eproof

\subsubsection{Remarks}
One could ask if statements (a) and (b) of Theorem
\ref{main_theorem} might
be equivalent in general. In more detail, is a $C^*$-algebra $A$ stable
if and only if $A$ admits no non-zero unital quotient and no non-zero
bounded trace? U.~Haagerup has proved that quasi-traces
on exact $C^*$-algebras are traces (cf.~\cite{Haa:quasi}), so
quasi-traces may be considered instead of traces. The answer to the
question is no. Based on a construction by J.~Villadsen
(cf.\ \cite{Vil:sr=n}), M.~Rørdam has produced an example of a simple,
nuclear $C^*$-algebra $B$ of stable rank one such that $M_2(B)$ is
stable but $B$ is not stable (cf.\ \cite{Ror:sns}). Suppose $A$ is a
$C^*$-algebra such that $M_2(A)$ is stable and $A$ is not stable. Then
$M_2(A)$ admits no bounded quasi-trace, and no quotient of $M_2(A)$ is
unital. This implies that $A$ admits no bounded quasi-trace and that
no quotient of $A$ is unital.

It would be interesting to determine the class of $C^*$-algebras for which
stability can be characterized by the non-existence of unital quotients and
bounded traces. It is noteworthy that for $C^*$-algebras in this class, the
extension problem for stability can be solved affirmatively, i.e., the
extension of two stable $C^*$-algebras is stable (cf.\ Corollary
\ref{thm3-1}). Also, for $C^*$-algebras in this class it follows that every
simple $C^*$-algebra is either stably finite or purely infinite (cf.\ 
\cite[Section 5]{HjeRor:stable}).

\mathversion{bold}
\subsection{Purely infinite $C^*$-algebras of locally finite graphs.} 
\mathversion{normal}

\begin{theorem}
\label{thm-1}
Let $E$ be a locally finite directed graph with no sinks. The following
six statements are equivalent:
\begin{itemize}
\item[(a)] $C^*(E)$ is purely infinite.
\item[(b)] $C^*(E)$ admits no non-zero trace.
\item[(c)] No quotient of $C^*(E)$ contains a two-sided closed ideal
  that is an $\AF$-algebra, nor does it contain a corner that is
  $*$-isomorphic to $M_n(C(\bT))$ for some $n\in\bN$.
\item[(d)] Every infinite path in $E$ admits a detour $\beta$ such that
  there are two or more cycles based at some vertex of ${\beta}^0$.
\item[(e)] Every vertex connects to a cycle with an exit in every subgraph
  $(F,\{e\mid r(e)\in F\})$, where $F$ is a complement of a hereditary and
  saturated subset of vertices in $E$.
\item[(f)] For every vertex $v\in E^0$ the projection $P_v\in C^*(E)$ is
properly infinite.
\end{itemize}
\end{theorem}

\bproof

{\bf (a) $\Rightarrow$ (b):} This follows from \cite[Proposition
5.1]{KR:purely}.  

{\bf (b) $\Rightarrow$ (c):} This follows by noting that an $\AF$-algebra
admits a densely defined trace and that a densely defined trace on a
two-sided closed ideal of a quotient of a $C^*$-algebra $A$ lift to a trace
on a two-sided ideal of $A$. Also, if $PAP$ is a corner in $A$ for some
projection $P\in A$ and if $\tau$ is a bounded everywhere defined trace on
$PAP$, then the map $\sum_{i=1}^n{x_iPaPy_i}\mapsto
\sum_{i=1}^n{\tau(PaPy_ix_iP)}$ extends $\tau$ to give a densely defined
trace on the two-sided ideal $\{\sum_{i=1}^n{x_iPaPy_i}\mid n\in\bN,
a,x_i,y_i\in A\}$ of A generated by $PAP$.

{\bf (c) $\Rightarrow$ (d):} First we show that no vertex of $E$ is
base for precisely one cycle. Let $v\in E^0$ and suppose to reach a
contradiction that precisely one cycle $\alpha$ is based at $v$. Set
\[
L(v) = \{w\in E^0\mid \mbox{$s(\gamma)=w$ and $r(\gamma)=v$ for some
  $\gamma\in E^*$}\},
\]
and let $H^0 = E^0\backslash L(v)$ denote the complement. Then $H^0$
is a saturated and hereditary subset of $E^0$ (see the proof of
Theorem \ref{main_theorem}) and if $F$ denotes the subgraph
($L(v),\{e\mid r(e)\in L(v)\}$) of $E$ then $C^*(F)$ is isomorphic to
a quotient of $C^*(E)$ by \cite[Theorem 6.6]{KPRR}. The cycle $\alpha$
based at $v$ is contained in $F$. Moreover, $\alpha$ has no exit in
$F$, hence it follows by Lemma \ref{torus} that $C^*(F)$ contains a
corner which is $*$-isomorphic to $M_n(C(\bT))$ for some $n\in\bN$.

We may assume that no vertex of $E$ is base for precisely one cycle.
Suppose now, to reach a contradiction, that there is an infinite path
$\alpha\in E^{\infty}$ such that for every detour $\beta$ of $\alpha$,
no cycle is based at any vertex of ${\beta}^0$. Set
\begin{equation}
K = \{w\in E^0\mid \mbox{$w\in\beta^0$ for a detour $\beta$ of $\alpha$}\}.
\end{equation} 
Note that $\a^0$ is contained in $K$. Set
\begin{equation}
H = \{w\in E^0\backslash K\mid \mbox{$s(\gamma)\in K$ and
$r(\gamma)=w$ for some $\gamma\in E^*$}\}.
\end{equation}
Let $\widetilde{H}$ be the smallest saturated subset containing the
hereditary completion of $H$. We show that $K$ and $\widetilde{H}$ are
disjoint sets. Note first that $K$ is disjoint to the hereditary
completion of $H$, and next that the smallest saturation of this set
does not contain any vertices of $K$, since for every vertex $w\in K$
there is an edge $e\in E^1$ such that $s(e)=w$, $r(e)\ne w$ and
$r(e)\in K$. Hence $K\cap \widetilde{H}=\emptyset$.
Let $K^{\prime}$ denote the smallest saturated subset containing $K$.
Then $K^{\prime}\cap \widetilde{H}=\emptyset$ because otherwise there
will be a path from $H$ to $K$ contradicting that $K\cap
\widetilde{H}=\emptyset$.  $K^{\prime}$ is hereditary in
$E^0\backslash\widetilde{H}$. It follows that $K^{\prime}$ is a
hereditary and saturated subset of $E^0\backslash\widetilde{H}$.
Moreover, no cycle is based at any vertex of $K^{\prime}$.

Now, at the $C^*$-algebra-level we obtain that if $F$ denotes the subgraph
$(E^0\backslash\widetilde{H},\{e\mid r(e)\in E^0\backslash\widetilde{H}\})$
then $C^*(F)$ is isomorphic to a quotient of $C^*(E)$. Let $I(K^{\prime})$
be the two-sided closed ideal of $C^*(F)$ induced by $K^{\prime}$ as in
\cite[Theorem 6.6]{KPRR}. Since no cycle is based at any vertex of
$K^{\prime}$ it follows that $I(K^{\prime})$ is an $\AF$-algebra by Theorem
\ref{af} or \cite[Theorem 2.4]{KPR}, which contradicts our assumption. This
shows that statement (c) implies (d).

{\bf (d) $\Rightarrow$ (e):} Let $H^0\subseteq E^0$ be a hereditary
and saturated subset and let $v\in E^0\backslash H^0$. Since $H^0$ is
saturated and $E$ has no sinks there is an infinite path $\a\in
E^{\infty}$ such that $v=s(\a)$ and $\a^0\subseteq E^0\backslash H^0$.
By assumption there is a detour $\beta$ of $\a$ such that two or more
cycles are based at some vertex of $\beta^0$. Since $H^0$ is
hereditary and $\a^0\subseteq E^0\backslash H^0$ it follows that
$\beta^0\subseteq E^0\backslash H^0$. Hence $v$ connects to a cycle
with an exit in the subgraph $(E^0\backslash H^0,\{e\mid r(e)\in
E^0\backslash H^0\})$.

{\bf (e) $\Rightarrow$ (f):} Let $v\in E^0$ be a vertex and let
$P_v\in C^*(E)$ be the associated projection. By \cite[Proposition
3.14]{KR:purely} $P_v$ is properly infinite when for every two-sided
closed ideal $I$ in $C^*(E)$, the projection $P_v+I$ is either zero or
infinite.

No vertex in $E^0$ is base for precisely one cycle, since otherwise
there is, as shown in ``(c) implies (d)'', a complement of a
hereditary and saturated subset of $E^0$ that contains a vertex which
is base for a cycle with no exit, and this contradicts the assumption.
Hence there is an isomorphism of the lattice of saturated hereditary
subsets of $E^0$ onto the lattice of ideals in $C^*(E)$ by
\cite[Theorem 6.6]{KPRR}.  Let $I$ be an ideal of $C^*(E)$ and let
$H^0$ be a saturated hereditary subset of $E^0$ such that $I=I(H^0)$.

If $v\in H^0$ then $P_v\in I$ and we are done.  Suppose $v\notin H^0$.
By assumption there is a path $\a$ such that $\alpha^0 \cap
H^0=\emptyset$, $s(\alpha)=v$ and there is a cycle $\gamma$ based at
$r(\a)$. Set $w=r(\a)$ and note that there is another cycle $\delta$
based at $w$ since no vertex in $E^0$ is base for precisely one cycle,
as shown above. Since $H^0$ is hereditary and $w\notin H^0$, the
subsets $\gamma^0$ and $\delta^0$ are both disjoint to $H^0$. The
projection $P_w$ does not belong to $I$ and $P_w$ is properly infinite
since
\[
S_{\gamma}^*S_{\gamma}=S_{\delta}^*S_{\delta}=P_w\geq 
S_{\gamma}S_{\gamma}^*+S_{\delta}S_{\delta}^*.
\]
Hence $P_v+I$ is (properly) infinite. Since $v$ connects to $w$ via
$\a$ with $\alpha^0 \cap H^0=\emptyset$, it follows that
$P_w+I\lesssim P_v+I$ and $P_v+I$ must be infinite.

{\bf (e) $\Rightarrow$ (a):} We aim to show that for every closed
two-sided ideal $I$ and every hereditary $C^*$-subalgebra $B$ in the
quotient by $I$, there is an infinite projection in $B$. Then $C^*(E)$
is purely infinite by \cite[Proposition 4.7]{KR:purely}.

Let $I$ be a two-sided closed ideal of $C^*(E)$. The assumption implies
that no vertex of $E$ is base for precisely one cycle (see ``(c) implies
(d)''). Hence $I=\mbox{Ideal}(H^0)$ for some hereditary saturated subset
$H^0\subseteq E^0$ (cf.\ \cite[Theorem 6.6]{KPRR}) and if $F$ denotes the
subgraph $(E^0\backslash H^0,\{e\mid r(e)\in E^0\backslash H^0\})$ then
$C^*(F)\cong C^*(E)/I$ by \cite[Theorem 6.6]{KPRR}. Notice that every
projection associated to a vertex of $F$ is properly infinite, since this
property is preserved by the quotient-map.

Let $B$ be a hereditary $C^*$-subalgebra of $C^*(F)$. For each
$n\in\bN$ the projections $\{S_{\a}S_{\a}^* \colon |\a|=n\}$ are
mutually orthogonal and span a $C^*$-subalgebra $A_n$ of $C^*(F)$.
Moreover, the projections are properly infinite since
$S_{\a}^*S_{\a}=P_{r(\a)}$. Set $D=(\cup A_n)^{-}$. By \cite[Lemma
3.5]{KPR} $C^*(F)$ contains a partial isometry $v$ such that $v^*v\in
D$ and $vv^*\in B$. Hence there is a (properly) infinite projection in
$A_n$ for some $n$ which is equivalent to a subprojection of $B$. This
shows that $C^*(E)$ is purely infinite and the proof is complete.
\eproof

\subsection{Examples}

\begin{description}

\item[Example (i):]
Consider the graph $B$ given by:
$$
B\colon\qquad
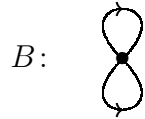
\xymatrix@=40pt{*-{\bullet}\pil@(ul,ur)[]\pil@(dl,dr)}$$
The Cuntz-Krieger $B$-family is given by the set
$\{S_e,S_f\colon\mbox{e,f are the edges of $B$}\}$ of partial
isometries subject to the relations
$$ S_e^*S_e=S_f^*S_f=S_eS_e^*+S_fS_f^* $$
Hence the generators of $C^*(B)$ satisfies the Cuntz-relations for $\cO_2$
and $C^*(B)\cong \cO_2$. $C^*(B)$ is purely infinite by Theorem \ref{thm-1}

%\item[Example (ii):]
%Consider the graph $D$ given by:
%$$
%D\colon\qquad
%\xymatrix@=40pt{
%  *-{\bullet}\pil@/^\boej/[r]\pil@/_\boej/[r]&
%  *-{\bullet}\pil@/^\boej/[r]\pil@/_\boej/[r]&
%  *-{\bullet}\pil@/^\boej/[r]\pil@/_\boej/[r]&
%  *-{\bullet}\pil@/^\boej/[r]\pil@/_\boej/[r]&
%  }\;\;\cdots$$
%\clearpage

%Represent $D$ by
%\begin{align*}
%D^0 &= \{0,1,2,\dots\},\qquad D^1 = \{e_k,f_k\colon k\ge 0\}\\
%s(e_k) &= k,\quad r(e_k)=k+1,\quad s(f_k)=k,\quad r(f_k)=k+1.
%\end{align*}
%Define $\tau_D\colon D^0\to\bR^+$ by $\tau_D(k)=2^{-k}$. Then $\tau_D$
%is a bounded graph-trace on $D$, and
%it follows from \cite{Hje:report} that $C^*(D)$ admits a
%bounded trace, hence is not stable.
%By \cite[Corollary 2.3]{KPR} $C^*(D)$ is the inductive limit
%given by
%\[
%M_3(\bC)\to M_7(\bC)\to M_{15}(\bC)\to M_{31}(\bC)\to\cdots\to C^*(D)
%\]
%where the connecting maps are embeddings of multiplicity two into the upper
%left-hand corner.

\item[Example (ii):]{\bf (cf. \cite{Kumjian:survey})}
Consider the graph $E$ given by:
$$\entrymodifiers={-{\bullet}}
E\colon\qquad
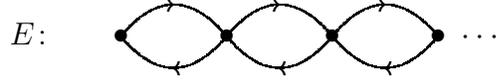
\xymatrix@=40pt{
  \pil@/^\boej/[r]&
  \pil@/^\boej/[r]\pil@/^\boej/[l]&
  \pil@/^\boej/[r]\pil@/^\boej/[l]&
%  \pil@/^\boej/[r]\pil@/^\boej/[l]&
  \pil@/^\boej/[l] }\;\;\cdots
$$ 

$C^*(E)$ is purely infinite by Theorem \ref{thm-1} and stable by Theorem
\ref{main_theorem}. Using the the Kirchberg-Phillips classification theorem
it can be shown that $C^*(E)\cong \cO_\infty\otimes\cK$. It is an open
question to the author whether this isomorphism can be given explicit.

%\item[Example (iv):]
%Consider the graph $F$ given by:
%$$
%F\colon\qquad
%\xymatrix@=40pt{
%  *-{\bullet}\pil@/^\boej/[r]\pil@/_\boej/[r]\pil[d]&
%  *-{\bullet}\pil@/^\boej/[r]\pil@/_\boej/[r]\pil[d]&
%  *-{\bullet}\pil@/^\boej/[r]\pil@/_\boej/[r]\pil[d]&
%  *-{\bullet}\pil@/^\boej/[r]\pil@/_\boej/[r]\pil[d]&
%  *-{\bullet}\pil[d]\ar@{}|{\textstyle\cdots}[r]&*{}\\
%  *-{\bullet}\pil@/^\boej/[r]&
%  *-{\bullet}\pil@/^\boej/[r]\pil@/^\boej/[l]&
%  *-{\bullet}\pil@/^\boej/[r]\pil@/^\boej/[l]&
%  *-{\bullet}\pil@/^\boej/[r]\pil@/^\boej/[l]&
%  *-{\bullet}\pil@/^\boej/[l]\ar@{}|{\textstyle\cdots}[r]&*{}
%  }$$
%By \cite[Theorem 6.6]{KPRR} the lower row of $F$ forms a hereditary and
%saturated subset which induces an ideal Morita equivalent to $C^*(E)$
%of Example (iii) and the quotient, induced by the upper row of $F$, is
%isomorphic to $C^*(D)$ of Example (ii). It
%follows that $C^*(F)$ admits a (bounded) trace and is not purely infinite
%by Theorem \ref{thm-1}. By \cite[Theorem 3.9]{KPR} 
%every hereditary $C^*$-subalgebra of $C^*(F)$ contains an infinite
%projection.

\end{description}

\mathversion{bold}
\section{Purely infinite $C^*$-algebras of groupoids.}
\mathversion{normal}
\label{groupoids}

In this section, $G$ is a second countable, locally compact and r-discrete
groupoid.  In \cite{A-Delaroche:purely} C.\ Anantharaman-Delaroche
introduces the notion of locally contracting r-discrete groupoids; $G$ is
{\em locally contracting} if for every non-empty open subset $U$ of the
unit space of $G$, there exist an open subset $V$ in $U$ and an open
bisection $S$ (i.e., $r$ and $d$ are one to one on $S$) with
$\overline{V}\subseteq d(S)$ and $\a_{S^{-1}}(\overline{V})$ strictly
contained in $V$, where $\a_{S}$ is the $G$-map associated with S, i.e.,
the homeomorphism of $r(S)$ onto $d(S)$ given by $\a_{S}(x)=d(xS)$ for all
$x\in r(S)$ (cf.\ \cite[Definition 2.1]{A-Delaroche:purely}).

In \cite[Proposition 2.4]{A-Delaroche:purely} it is shown that every
non-zero hereditary $C^*$-subalgebra of the reduced $C^*$-algebra of a
locally contracting groupoid $G$ contains an infinite projection, when the
elements of the unit space of $G$ for which the isotropy-group is trivial,
are dense.  When the latter condition holds for every closed invariant
subset of the unit space, $G$ is said to be {\em essentially principal}
(cf.\ \cite[Definition 4.3]{Ren:groupoid}).

\begin{proposition}
\label{pi_groupoid}
Let $G$ be an r-discrete and essentially principal groupoid. Assume for
every closed invariant subset $F$ of the unit space that the reduced
subgroupoid $G_F$ is locally contracting. Then $C_r^*(G)$ is purely
infinite.
\end{proposition}

\bproof
$C_r^*(G)$ is purely infinite if every non-zero hereditary $C^*$-subalgebra
in every quotient contains an infinite projection (cf.\ \cite[Proposition
4.7]{KR:purely}). Let $I$ be a two-sided closed ideal of $C_r^*(G)$. By
\cite[Proposition 4.6]{Ren:groupoid} $I$ is the closure of the
$*$-algebra of continuous functions with compact support that vanishes
outside $U$ for some open invariant subset $U\subseteq G^0$, and 
\[
C_r^*(G_F)\cong C_r^*(G)/I,
\]
where $F=G^0\backslash U$ is a closed invariant subset of $G^0$. By
assumption the elements of $G_F^0$ (=$F$) with trivial isotropy
are dense and $G_F$ is locally contracting. It follows from
\cite[Proposition 2.4]{A-Delaroche:purely} that every non-zero hereditary
$C^*$-subalgebra of $C_r^*(G_F)$ contains an infinite projection, hence
$C_r^*(G)$ is purely infinite.
\eproof

\subsection{Singly generated dynamical systems}
\label{sgds}

Following J.\ Renault \cite{Ren:cuntz-like} (cf.\ \cite{Rome-treaty}) a singly
generated dynamical system (SGDS), denoted $(X,T)$, consists in our setting
of a locally compact, second countable Hausdorff space $X$ and a local
homeomorphism $T$ from an open subset dom($T$) of $X$ onto an open subset
ran($T$) of $X$. The semi-direct product groupoid $G(X,T)$ attached to
$(X,T)$ is defined as follows. Let 
\[
G(X,T)=\{(x,m-n,y)\mid \mbox{$(x,y)\in$ dom($T^m$)$\times$dom($T^n$),
  $T^m(x)=T^n(y)$ for 
  $n,m\in \bN$}\}
\]
with groupoid structure:
\[
(x,k,y)(y,l,z)=(x,k+l,z),\qquad (x,k,y)^{-1}=(y,-k,x)
\]
and topology defined by the basic open sets
\[
{\cal U}(U,m,n,V)=\{(x,m-n,y)\mid \mbox{$(x,y)\in U\times V$ such that
  $T^m(x)=T^n(y)$}\}
\]
where $U$ and $V$ are open subsets of dom($T^m$) resp.\ dom($T^n$) on which
$T^m$ resp.\ $T^n$ are injective. The range and source maps, 
$r(x,m-n,y)=x$ and $d(x,m-n,y)=y$, are local homeomorphisms and in this
case the groupoid is called ``\'etale''. We remark that in case there is no
non-empty open set on which $T^m$ and $T^n$ agree for all distinct numbers
$m,n\in \bN$ the groupoid $G(X,T)$ is essentially free and isomorphic to
the groupoid of germs of the pseudogroup
generated by the restrictions $T_{|U}$, where $U$ is an open subset of
$X$ (cf.\ \cite{Rome-treaty}).

The groupoid $G(X,T)$ of an SGDS $(X,T)$ is a locally compact \'etale Hausdorff
groupoid that is also amenable by
\cite[Proposition 2.4]{Ren:cuntz-like} . Hence the full and reduced
$C^*$-algebra $C^*(X,T)$ of the convolution algebra of $G(X,T)$ constructed in
\cite{Ren:groupoid} coincide and $C^*(X,T)$ is nuclear (see
\cite{RenA-Delaroche:amenable}). We can express in 
case of an SGDS $(X,T)$ the conditions of Proposition \ref{pi_groupoid}
ensuring that the $C^*$-algebra is purely infinite.

\begin{corollary}
\label{pi_sgds}
Let $(X,T)$ be an SGDS. Assume for every closed invariant subset
$F\subseteq X$ and for every open subset $U\subseteq F$ that
\begin{itemize}
\item[(i)] There is an element $u\in U$ such that $T^n(u)\neq T^m(u)$ for
  all distinct numbers $n,m\in\bN$.
\item[(ii)] There is an open set $V\subseteq U$ and distinct numbers
  $n,m\in\bN$ such that $T^n(\overline{V})$ is strictly contained in $T^m(V)$.
\end{itemize}
Then $C^*(X,T)$ is purely infinite.
\end{corollary}

\bproof
Let $U$ and $F$ be as in the corollary. It follows from condition (i)
that every open subset of $F$ contains an element with trivial isotropy, hence
$G(X,T)$ is essentially principal. 

Let $V\subseteq U$ and $n,m\in\bN$ be as in condition (ii). We may assume
that $T^n$ and $T^m$ are both injective on $U$ and also that $\overline{V}$
is contained in some open subset $W\subseteq U$. One defines a bisection in
$G(X,T)_F$ by letting 
\[
S=\{(x,m-n,y)\mid\mbox{$T^m(x)=T^n(y)$ for $x,y\in W$}\}.
\]
Then $G(X,T)_F$ is locally contracting since $\a_{S^-1}(\overline{V})$
equals $T^{n-m}(\overline{V})$ and hence is strictly contained in $V$ by
assumption. It follows from Proposition \ref{pi_groupoid} that $C^*(X,T)$
is purely infinite. 
\eproof

\begin{remark}
  For the locally finite graphs $E$ (without sinks) of section
  \ref{purely}, the construction in \cite{KPRR} of $C^*(E)$ reveals that
  this $C^*$-algebra comes from an SGDS $(X_{A_E},T_{A_E})$, where
  $T_{A_E}$ is the one-sided shift on the infinite path space $X_{A_E}$
  associated to the edge-matrix $A_E$ of the graph $E$. Corollary
  \ref{pi_sgds} gives an alternative proof of ``(e) implies (a)'' of
  Theorem \ref{thm-1}.
\end{remark}

\mathversion{bold}
\subsection{$C^*$-algebras of infinite matrices.}
\mathversion{normal}
\label{infinite_matrices}

In the following $A$ will denote an $I\times I$-matrix with entries in
$\{0,1\}$ for an at most countable set $I$. The associated graph $E_A$ to
$A$ is the directed graph with vertex set $I$ and edge set the subset of
$I\times I$ for which $A(i,j)=1$. Paths in $E_A$ consist of sequences of
vertices $(i_k)\subseteq I^{\bN}$ such that $A(i_k,i_{k+1})=1$. A cycle is a
path $(i_1,\ldots,i_n)$ such that $A(i_n,i_1)=1$, and a cycle has an exit if
there is a vertex $j\in I$ such that $A(i_k,j)=1$ and $j\neq i_{k+1}$ for
some vertex $i_k$ of the cycle (where $i_{n+1}=i_1$).

Assuming that $A$ has no zero-rows, a generalized Cuntz-Krieger $C^*$-algebra
${\cO}_A$ is constructed in \cite{ExelLaca:infinite} as follows: 
Suppose that $\{S_i\mid i\in I\}$ is a set of partial isometries indexed by
the vertex set $I$ such the the range projections are mutually orthogonal and 
\begin{itemize}
\item[(i)] $S_i^*S_i$ commutes with $S_j^*S_j$ for all $i,j$,
\item[(ii)] $S_i^*S_j=0$ if $i\neq j$,
\item[(iii)] $(S_i^*S_i)S_j=A(i,j)S_j$ for all $i,j$,
\item[(vi)] For each pair of finite subsets $X,Y\subseteq I$ such that
\[
A(X,Y,j):=\prod_{x\in X}A(x,j)\prod_{y\in Y}(1-A(y,j))
\]
vanishes for all but finitely many $j$'s, the equation
\[
\prod_{x\in X}S_x^*S_x\prod_{y\in Y}(1-S_y^*S_y)=\sum_{j\in I}A(X,Y,j)S_jS_j^*
\]
holds.
\end{itemize} 
The universal $C^*$-algebra generated by such a set $\{S_i\mid i\in I\}$ is
denoted by ${\cO}_A$.

Note that for every pair $i,j\in I$, the element $S_iS_j$ is non-zero if
and only if the entry $A(i,j)$ is non-zero, i.e., there is an edge from
vertex $i$ to vertex $j$. If $\a=(i_1,\ldots,i_n)$ is a path, we define a
partial isometry by $S_{\a}:=S_{i_1}\cdots S_{i_n}$. Using relation (iii)
stated above it follows for finite paths $\a$ and $\b$ that
\begin{equation} \label{eq: path}
S_{\b}^*S_{\a} = \left\{ \begin{array}{ll}
    S_{\b'}^* & \textrm{if $\b=\a\b'$}\\
    S_{\b}^*S_{\b} & \textrm{if $\b=\a$}\\
    S_{\a'} & \textrm{if $\a=\b\a'$}\\
    0 & \textrm{otherwise}
     \end{array} \right.
\end{equation}
(c.f.\ \cite[Lemma 1.1]{KPR}). Set
\begin{eqnarray*}
M=\{S_{\a}(\prod_{i\in X}S_i^*S_i)S_{\b}^*
\mid \a ,\b\in E_A^*,\ X\subseteq I, \mbox{finite}\},
\end{eqnarray*}
where $S_{\emptyset}=1$ and note that $S_i=S_iS_i^*S_i\in M$ for all $i\in
I$. A computation using  
(\ref{eq: path}) and that initial projections commute, then reveals that 
Span$\{M\}$ is an algebra. Also, Span$\{M\}$ is closed under taking adjoints 
and ${\cO}_A$ is indeed the norm-closure of Span$\{M\}$.

The theorem that follows below, characterizes $AF$-algebras among
the $C^*$-algebras ${\cO}_A$ and is a generalization of 
\cite[Theorem 2.4]{KPR}. First we state a lemma.

\begin{lemma} \label{torus}
Let $A$ be an $I\times I$-matrix with no zero rows and let $E_A$ be the 
associated graph. If $E_A$ has a cycle of period $n$ with no exits, then
there is a projection $P\in {\cO}_A$ such that the corner $P{\cO}_AP$ is 
$*$-isomorphic to $M_n(C(\bT))$.
\end{lemma}

\bproof
Let $\a=(i_1,\ldots,i_n)$ be a cycle of period $n$ with no exits, i.e.,
$i_{n+1}=i_1$, $A(i_n,i_1)=1$ and if $A(i_k,j)=1$ then $j=i_{k+1}$ for
$j\in I$ and $k=1,\ldots,n$. The associated partial
isometries $S_{i_1},\ldots,S_{i_n}$ satisfy that
$S_{i_k}^*S_{i_k}=S_{i_{k+1}}S_{i_{k+1}}^*$, hence
the initial projections $S_{i_k}^*S_{i_k}$ are mutually orthogonal and the
partial isometry $S_{\a}:=S_{i_1}\cdots S_{i_n}$ satisfies that
$S_{\a}^*S_{\a}=S_{\a}S_{\a}^*$. Since for every set
$\{c_{i_1},\ldots,c_{i_n}\}$ of scalars of norm 1, the universal
$C^*$-algebra generated by $\{c_{i_1}S_{i_1},\ldots,c_{i_n}S_{i_n}\}$ is
$*$-isomorphic to $C^*(S_{i_1},\ldots,S_{i_n})$, it follows that $\bT$ is
contained in the spectrum of $S_{\a}$. Hence $S_{\a}$ is a partial unitary
with full spectrum, whence $C^*(S_{\a})\cong C(\bT)$. There is a
$*$-isomorphism of $C^*(S_{i_1},\ldots,S_{i_n})$ onto $C^*(S_{\a})\otimes M_n$
given by
\begin{equation*} 
S_{i_k}\mapsto \left\{ \begin{array}{ll}
    1\otimes e_{k,k+1} & \textrm{if $k\in\{1,\ldots,n-1\}$,}\\
    S_{\a}\otimes e_{n,1} & \textrm{$k=n$,}
     \end{array} \right. 
\end{equation*}
where $1$ is the unit of $C^*(S_{\a})$ and $e_{j,k}$ are matrixunits (cf.\
\cite [Remark 3.9]{CK:algebras}). Hence $C^*(S_{i_1},\ldots,S_{i_n})\cong
C(\bT)\otimes M_n$.

Set $P=\sum_{k=1}^nS_{i_k}S_{i_k}^*$ and note that $PS_j=S_j=S_jP$ if 
$j\in\{i_1,\ldots,i_n\}$ and $PS_j=0=S_jP$ otherwise. Also,
\begin{equation*} 
P(\prod_{i\in X}S_i^*S_i)P = \left\{ \begin{array}{ll}
      0 & \textrm{if $A(i,i_k)=0$, $i\in X$, $i_k\in\{i_1,\ldots,i_n\}$}\\ 
      \sum_{k\in J}^nS_{i_k}S_{i_k}^* &
      \textrm{otherwise. ($J\subseteq\{1,\ldots,n\}$).} 
     \end{array} \right.
\end{equation*}
Hence for every finite subset $X\subseteq I$ and every pair of finite paths
$\a,\b\in E_A^*$ (including the empty path) it follows that
\[
PS_{\a}(\prod_{i\in X}S_i^*S_i)S_{\b}^*P\in C^*(S_{i_1},\ldots,S_{i_n}).
\]
This shows that the corner $P{\cO}_AP$ is contained in
$C^*(S_{i_1},\ldots,S_{i_n})$. Using once more that $PS_j=S_j=S_jP$ if 
$j\in\{i_1,\ldots,i_n\}$,  $P{\cO}_AP$ equals $C^*(S_{i_1},\ldots,S_{i_n})$
and is thus $*$-isomorphic to $M_n(C(\bT))$.
\eproof

\begin{remark}
  In case of row-finite graphs $E$ with no sinks, the argument given in the
  proof of Lemma \ref{torus} applies such that if $E$ has a cycle of period
  $n$ with no exit, then $C^*(E)$ contains a projection $P$ such that the
  corner $PC^*(E)P$ is $*$-isomorphic to $M_n(C(\bT))$.
\end{remark}

\begin{theorem}
\label{af}
Let $A$ be an $I\times I$-matrix with no zero rows and let $E_A$ be the 
associated graph. Then ${\cO}_A$ is an $\AF$-algebra if and only if $E_A$ 
has no cycles.
\end{theorem}

\bproof Suppose that $E_A$ has no cycles.
We show that ${\cO}_A$ has the local finiteness property, i.e., that every 
finite set of elements in ${\cO}_A$ can be approximated by finitely many 
elements in a finite dimensional $C^*$-subalgebra of ${\cO}_A$. Let 
$F\subseteq I$ be a finite subset and set $B=C^*(\{S_i\mid i\in F\})$.
Then $B$ is a $C^*$-subalgebra of ${\cO}_A$ and we show that $B$ is finite 
dimensional.   
Set
\begin{eqnarray*}
M_F=\{S_{\a}(\prod_{i\in F_0}S_i^*S_i)S_{\b}^*
\mid \a ,\b\in F^*,F_0\subseteq F\},
\end{eqnarray*}
By the argument given above it follows that Span$\{M_F\}$ is a $*$-algebra 
containing $S_i$ for all $i\in F$. Since $F$ is finite and the graph $E_A$ 
has no cycles, there can only be finitely many paths with vertices in $F$, 
hence Span$\{M\}$ is finite dimensional and equals $B$. This shows that
${\cO}_A$ is an $\AF$-algebra.

For the converse statement, suppose that $\a$ is a cycle in $E_A$ of period
$n$. If $\a$ has an exit then ${\cO}_A$ contains an infinite projection by
the argument given in the proof of \cite[Theorem 2.4]{KPR}, hence ${\cO}_A$
can not be an $\AF$-algebra. If $\a$ has no exits then there is a
projection $P\in {\cO}_A$ such that the corner $P{\cO}_AP$ is
$*$-isomorphic to $M_n(C(\bT))$ by Lemma \ref{torus}. A corner of an
$\AF$-algebra is an $\AF$-algebra, but $M_n(C(\bT))$ is
not $\AF$ (since $\bT$ is not totally disconnected) and hence ${\cO}_A$ is not
an $\AF$-algebra in this case either. This completes the proof.  
\eproof

\begin{remark} \label{af_remark}
  Theorem \ref{af} is inspired by \cite[Theorem 2.4]{KPR} and is easily
  translated to give a new proof in case of row-finite graphs $E$ with no
  sinks such that $C^*(E)$ is an $AF$-algebra if and only if $E$ has no
  cycles (cf.\ \cite[Theorem 2.4]{KPR}).
\end{remark}

\mathversion{bold}
\subsection{Pure infiniteness.}
\mathversion{normal}
\label{pure_infiniteness}

Theorem \ref{markov_pi} of this section is a characterization of pure
infiniteness for the $C^*$-algebras ${\cO}_A$ of section
\ref{infinite_matrices}.

In \cite[Section 4]{Ren:cuntz-like} J.~Renault constructs a groupoid model for 
${\cO}_A$. More precisely, there is an SGDS $(X_A,T_A)$, called a Markov
shift (cf.\ \cite[Definition 4.1]{Ren:cuntz-like}), such that
$C^*(X_A,T_A)={\cO}_A$. In \cite{ExelLaca:infinite} one distinguishes 
between ${\cO}_A$ and $\tilde{{\cO}}_A$, where the latter is the unital
$C^*$-algebra generated by $\{S_i\mid i\in I\}$ and $1$. Either ${\cO}_A
=\tilde{{\cO}}_A$ or $\tilde{{\cO}}_A$ is ${\cO}_A$ with a unit adjoined. 
Both $C^*$-algebras come from a Markov shift as shown in
\cite[Proposition 4.8]{Ren:cuntz-like}; ${\cO}_A$ comes from
$(X_A,T_A)$, where either one point is removed from $X_A$ or no point is
removed. 

We refer to \cite[Section 4]{Ren:cuntz-like} for the construction of $X_A$,
the compact and totally disconnected Hausdorff space consisting of {\em
  terminal paths} associated to an infinite matrix $A$, where a terminal
path is either an infinite path $(i_0,i_1,\ldots )$ or a path with finite
exit time (cf.\ \cite[Proposition 4.1]{Ren:cuntz-like}).  The partial
homeomorphism $T_A$ is the one-sided shift with domain a dense open subset
$U\subseteq X_A$ and range an open subset of $X_A$. The domain $U$ admits a
partition $\{U_i\mid i\in I\}$ of pairwise disjoint compact open subsets
$U_i$. For each $i\in I$, $U_i$ is the set of all terminal paths of lenght
$\ge 1$ starting with $i$.

%For a terminal path $x\in U$ there is an orbit $\{x,T_Ax,T_A^2x,\ldots\}$. The
%exit time of $x$ is determined by the map $\tau:X_A\to\{0,1,\ldots
%,\infty\}$ such that $\tau(x)$ is the smallest number, possibly $\infty$,
%for which $T_A^{\tau(x)}x\notin U$. If $\tau(x)=\infty$ then
%$x=(i_0,i_1,\ldots )$ is just an infinite path. If $\tau(x)=n$ then 
%$x=(i_0,i_1,\ldots,i_n;I(i_n))$.

\begin{remark} \label{isolated_periodic}
  A terminal path in $X_A$ is called {\em periodic} if it, after some
  iterations by $T_A$, repeats indefinitely a cycle. Moreover, if the cycle
  has no exit, then the path is called {\em isolated periodic}. From Lemma
  \ref{torus} it follows for a Markov shift $(X_A,T_A)$ of an $I\times
  I$-matrix $A$ with no zero rows that if $X_A$ has an isolated periodic
  path, then there is a corner in $C^*(X_A,T_A)$ which is $*$-isomorphic to
  $M_n(C(\bT))$. We express Theorem \ref{af} in terms of Markov shifts:
\end{remark}

\begin{corollary}
\label{markov_af}
Let $A$ be an $I\times I$-matrix with no zero rows and let $(X_A,T_A)$ be the 
associated Markov shift. Then $C^*(X_A,T_A)$ is an $\AF$-algebra if and
only if $X_A$ has no periodic paths.
\end{corollary}

\begin{theorem}
\label{markov_pi}
Let $A$ be an $I\times I$-matrix with no zero rows and let $(X_A,T_A)$ be the 
associated Markov shift. The following five statements are equivalent.
\begin{itemize}
\item[(a)] $C^*(X_A,T_A)$ is purely infinite.
\item[(b)] $C^*(X_A,T_A)$ admits no non-zero trace.
\item[(c)] No quotient of $C^*(X_A,T_A)$ contains a two-sided closed ideal
  that is an $\AF$-algebra, nor does it contain a corner that is
  $*$-isomorphic to $M_n(C(\bT))$ for some $n\in\bN$.
\item[(d)] For every closed invariant subset $F\subseteq X_A$, every vertex
  $k\in I$ of every terminal path in $F$ connects to a cycle with an exit,
  i.e., there are subpaths $(i_1,\ldots,i_n)$ and $\beta =(j_1,\ldots,j_m)$
  of paths in $F$ such that $k=i_1$, $i_n=j_1=j_m$ for some $m,n\in\bN$ and
  $\beta$ has an exit in $F$.
  \item[(e)] For every closed invariant subset $F\subseteq X_A$ and for every
  open subset $V\subseteq F$ the following conditions hold;
\begin{itemize}
\item[(i)] There is a terminal path $x\in V$ such that $T_A^n(x)\neq
  T_A^m(x)$ for all distinct numbers $n,m\in\bN$.
\item[(ii)] There is an open set
  $W\subseteq V$ and distinct numbers $n,m\in\bN$ such that
  $T_A^n(\overline{W})$ is strictly contained in $T_A^m(W)$. 
\end{itemize}
\end{itemize}
\end{theorem}

\bproof

{\bf (a) $\Rightarrow$ (b):} This follows from \cite[Proposition
5.1]{KR:purely}.   

{\bf (b) $\Rightarrow$ (c):} See (b) implies (c) of Theorem \ref{thm-1}.

{\bf (c) $\Rightarrow$ (d):} Let $F\subseteq X_A$ be a closed invariant
subset and let $G(X_A,T_A)_F$ denote the reduced subgroupoid of
$G(X_A,T_A)$. Suppose to reach a contradiction that there is a vertex of a
path in $F$ which connects to a cycle with no exit (in sense of (d)).
Hence $F$ contains an isolated periodic path, whence the quotient
$C^*(G(X_A,T_A)_F)$ of $C^*(X_A,T_A)$ contains a corner which is
$*$-isomorphic to $M_n(C(\bT))$ for some $n\in\bN$ (cf.\ Lemma \ref{torus}
and Remark \ref{isolated_periodic}). This contradicts statement (c) and we
may assume that every cycle has an exit.

Suppose now that there is a vertex $k\in I$ of some path in $F$ such that
$k$ does not connect to a cycle. Let $V_k$ denote the open subset of $F$
consisting of all terminal paths in $F$ starting with $k$, and let $O(V_k)$
denote the union of all orbits $T_A^n(V_k)$ of $V_k$. Set
\[
{\tilde O}(V_k)=\{x\in F\mid \mbox{$T_A^m(x)\in O(V_k)$ for some $m\in\bN$}\}
\]
Then $\tilde{O}(V_k)$ is an open invariant subset of $F$ without periodic
paths. By \cite[Proposition 4.5]{Ren:groupoid}, $\tilde{O}(V_k)$ induce a
two-sided closed ideal of the quotient $C^*(G(X_A,T_A)_F)$, and this ideal
is an $\AF$-algebra by Corollary \ref{markov_af}. This contradicts the
assumption and statement (d) holds.

{\bf (d) $\Rightarrow$ (e):} Let $F\subseteq X_A$ be a closed invariant
subset and let $V$ be an open subset of $F$. By assumption every cycle has
an exit within $F$. There is a finite path $(i_1,\ldots,i_n)$ such that the
cylinderset $Z(i_1,\ldots,i_n)\subseteq V$ of all terminal paths starting
with $(i_1,\ldots,i_n)$ is an open subset of $V$. Since every cycle has an
exit, there is a path in $Z(i_1,\ldots,i_n)$ which is not periodic, hence
condition (i) of statement (e) holds.  Since every vertex of every path in
$F$ connects to a cycle with an exit, the argument given in the proof of
\cite[Proposition 4.10]{Ren:cuntz-like} shows that condition (ii) holds.

{\bf (e) $\Rightarrow$ (a):} This follows from Corollary \ref{pi_sgds} and
the proof is complete.

\eproof

\providecommand{\bysame}{\leavevmode\hbox to3em{\hrulefill}\thinspace}

\vspace{.5cm}
\noindent{\sc Department of Mathematics and Computer Science, Odense
University,
Campusvej 55, DK-5230 Odense M, Denmark}\\

\noindent{\sl E-mail address:} {\tt jacobh@imada.sdu.dk}\\

\end{document}